\documentstyle[11pt]{article}
\begin{document}
\newcommand{\p}{\parallel }
\makeatletter \makeatother
\newtheorem{th}{Theorem}[section]
\newtheorem{lem}{Lemma}[section]
\newtheorem{de}{Definition}[section]
\newtheorem{rem}{Remark}[section]
\newtheorem{cor}{Corollary}[section]
\renewcommand{\theequation}{\thesection.\arabic {equation}}

\title{{\bf Adiabatic Limits, Vanishing Theorems and the Noncommutative
Residue}
\author{ Kefeng Liu, Yong Wang}}

\date{}
\maketitle

\begin{abstract} In this paper, we compute the adiabatic limit of the scalar curvature and prove
several vanishing theorems, we also derive a Kastler-Kalau-Walze
type theorem for the noncommutative residue in the case of
foliations.
 \\

\indent{\bf Keywords:}\quad
 Foliations; adiabatic limits; vanishing theorems; noncommutative residue\\
\indent {\bf 2000 MSC:} 53C27, 51H25, 46L87\\

\end{abstract}

\section{Introduction}

\quad Let $(M,F)$ be a compact foliated manifold. If $F$ is spin and
there is a metric $g^F$ on $F$ such that the scalar curvature of
$g^F$ is positive, then $\widehat{A}$-genus of $M$ vanishes as
proved in [Co]. A. Connes proved this theorem by a highly
noncommutative method. In [LZ], Liu and Zhang used adiabatic limits
to study foliated manifolds. They constructed the sub-Dirac operator
associated to foliations with spin leave and presented a direct
geometric proof of the vanishing theorem of Connes for almost
Riemannian foliations by taking the adiabatic limit of the square of
the sub-Dirac operator. In [LMZ], the author used the sub-Dirac
operator to prove rigidity theorems for foliations. In [Ko2], [LK],
adiabatic limits were used to study the spectral geometry for
Riemannian foliations. In [Ru], Rumin studied the adiabatic limits
of some geometric objects for contact
manifolds.\\
\indent Here are some motivations of this paper.\\
\indent 1) For contact manifolds and in general for a manifold $M$
with splitting tangent bundle $TM=F\oplus F^{\bot}$ and $F$ is not
integrable, we hope to get a
Connes type vanishing theorem by taking adiabatic limits.\\
\indent 2) For complex foliations, by using the
Bismut-Kodiara-Nakano formula in [Bi] and taking adiabatic limits,
we hope to get a
vanishing theorem for foliation.\\
\indent 3) The Kastler-Kalau-Walze theorem says that the
noncommutative residue (see [Wo], [FGV]) of the $-{\rm dim}M+2$
power of the Dirac operator for even-dimensional spin manifolds $M$
is proportional to the Einstein-Hilbert action. For a foliation
$(M,F)$ with spin leave, by using the sub-Dirac operator in [LZ] and
considering the adiabatic limit of the noncommutative residue of the
$-{\rm dim}M+2$ power of the sub-Dirac operator, we hope to derive
a Kastler-Kalau-Walze type theorem for foliations.\\
 \indent In this paper, We define $\Phi(\omega),~A,~B,~ \Psi$ by
 (2.12), (2.20), (3.24), (4.33) in the following respectively. We will prove the following theorems.\\

 \noindent{\bf Theorem I}~~{\it Let $(M,F,g^F)$ be a compact and
transversally oriented foliation with spin leave, if
 $A(\omega,\phi_i (F^{\bot}))>0$ for any $\phi_i(F^{\bot})$ appeared
in (2.22), then
$<\widehat{A}(TM),[M]>=0.$}\\

 \noindent{\bf Theorem II}~~{\it Let $M$ be a compact
oriented manifold, if 1) $TM=F\oplus F^{\bot}$ and $F$ may not be
integrable and oriented spin, 2) $B>0$, then
$<\widehat{A}(TM),[M]>=0.$}\\

\noindent{\bf Theorem III}~~{\it Let $(M,F)$ be a compact complex
 foliation. If $\Psi>0,$ then the Euler number ${\rm Eul}(\xi)$ of the holomorphic
 bundle $\xi$ vanishes.}\\

\noindent{\bf Theorem IV}~~{\it Let $(M^n,F)$ be a compact
even-dimensional oriented
 foliation with spin leave and codimension $q$, and $D_F$ be the sub-Dirac operator, then
${\rm lim}_{\varepsilon\rightarrow 0}\varepsilon^{\frac{q}{2}}{\rm
Res}(D^{-n+2}_{F,\varepsilon})$ is proportional to $\int_M[k^F+\Phi(\omega)]d{\rm vol}_{g}.$}\\

\indent This paper is organized as follows: In Section 2, we compute
the adiabatic limit of the scalar curvature explicitly by following
the method in [LZ], then by using this result we give a new
vanishing theorem for general foliations with spin leave. In Section
3, for a manifold $M$ with splitting tangent bundle $TM=F\oplus
F^{\bot}$ and $F$ may not be integrable, we also compute the
adiabatic limit of the scalar curvature similarly, and we note that
some extra singular terms $O(\frac{1}{\varepsilon})$ will appear. We
consider the adiabatic limit of $\varepsilon D^2_F$ and derive a
vanishing theorem. In Section 4, by the Bismut-Kodiara-Nakano
formula and taking adiabatic limits, a vanishing theorem can be
obtained for complex foliations. In Section 5, a formula similar to
the Kastler-Kalau-Walze
theorem for foliations with spin leave is given.\\

\section{Vanishing theorem for foliations }

 \quad First we recall the basic setup and some facts in [LZ] (for details, see
 [LZ]).\\
  \indent Let $(M,F)$ be a foliation, that is, $F$ is an integrable
 sub-bundle of the tangent bundle $TM$. Take a metric $g^{TM}$ on
 $TM$ as in [LZ], then
 $$TM=F\oplus F^{\bot};~g^{TM}=g^F\oplus g^{F^{\bot}},\eqno(2.1)$$
where $F^{\bot}$ is the orthogonal complement of $F$ in $TM$ with
respect to $g^{TM}$ and $g^F$ $({\rm resp.}~ g^{F^{\bot}})$ the
metric on $F$ $({\rm resp.}~ F^{\bot}).$ Let $p$, $p^{\bot}$ be the
orthogonal projection from $TM$ to $F,~F^{\bot}$. Let $\nabla^{TM}$
be the Levi-Civita connection of $g^{TM}$ and $\nabla^F$
(resp.$\nabla ^{F^{\bot}}$) be the restriction of $\nabla^{TM}$ to
$F$ (resp. $F^{\bot}$). For any $\varepsilon>0$. let
$g^{TM,\varepsilon}$ be the metric
$$g^{TM,\varepsilon}=g^F\oplus\frac{1}{\varepsilon}g^{F^{\bot}}.\eqno(2.2)$$
Let $\nabla^{TM,\varepsilon}$ be the Levi-Civita connection of
$g^{TM,\varepsilon}$ and $\nabla^{F,\varepsilon}$ (resp.$\nabla
^{F^{\bot},\varepsilon}$) be the restriction of
$\nabla^{TM,\varepsilon}$ to $F$ (resp. $F^{\bot}$). Then we have
formulas (1.5)-(1.8) in [LZ]. Let $\dot{\nabla}$ be the Bott
connection on $F^{\bot}$ and $\dot{\nabla}^{\star}$ be the dual
connection of $\dot{\nabla}$ and
$$\omega:=\dot{\nabla}^{\star}-\dot{\nabla};~
\widehat{\nabla}=\frac{\dot{\nabla}+\dot{\nabla}^{\star}}{2}.
\eqno(2.3)$$
Let $k^{TM,\varepsilon}$ be the scalar curvature associated to
$\nabla^{TM,\varepsilon}.$ In the following, we will compute the
adiabatic limit ${\rm lim}_{\varepsilon\rightarrow
0}(k^{TM,\varepsilon}).$

  Let $\{f_i\}^p_{i=1},\{h_s\}_{s=1}^q$ be an orthonormal
basis of $g^{TM}=g^F\oplus g^{F^{\bot}}.$ By (2.32) in [LZ], we get
    $${\rm lim}_{\varepsilon\rightarrow
 0}\langle R^{TM,\varepsilon}(f_i,f_j)f_i,f_j\rangle=\langle
 R^{F}(f_i,f_j)f_i,f_j\rangle.\eqno(2.4)$$
 By (2.35) in [LZ], for $X\in \Gamma
(F),~U, V\in \Gamma (F^{\bot}),$
 $$A:={\rm lim}_{\varepsilon\rightarrow
 0}p^{\bot}\nabla^{TM,\varepsilon}_UX=\frac{1}{2}\sum^q_{s=1}\omega(X)(U,h_s)h_s.\eqno(2.5)$$
By (2.5) and (2.33), (1.13) in [LZ], we get
$${\rm lim}_{\varepsilon\rightarrow
 0}\langle
 R^{TM,\varepsilon}(X,U)X,U\rangle=-\frac{1}{2}\omega(p\nabla^{TM}_XX)(U,U)-\frac{1}{2}\omega(X)(U,p^{\bot}[X,U])$$
$$+\langle[X,A],U\rangle+\frac{1}{2}\omega(X)(U,A).\eqno(2.6)$$
 By (1.7), (1.8), (1.6) and (1.13) in [LZ],
we have
     $$\langle\nabla^{TM,\varepsilon}_V\nabla^{TM,\varepsilon}_UU,V\rangle=
\langle\nabla^{TM,\varepsilon}_V(p\nabla^{TM,\varepsilon}_UU),V\rangle
+\langle
p^{\bot}\nabla^{TM,\varepsilon}_Vp^{\bot}\nabla^{TM,\varepsilon}_UU,V\rangle$$
$$=-\frac{1}{2}\langle
p\nabla^{TM,\varepsilon}_UU,2\nabla^{TM}_VV\rangle+\frac{\varepsilon}{2}\langle
p\nabla^{TM,\varepsilon}_UU,[V,V]\rangle+\langle
\nabla_V^{F^{\bot}}\nabla_U^{F^{\bot}},V\rangle$$
$$=-\langle
p\nabla^{TM,\varepsilon}_UU,\nabla^{TM}_VV\rangle+\langle
\nabla_V^{F^{\bot}}\nabla^{F^{\bot}}_U,V\rangle$$
$$=-\langle\nabla^{TM}_UU,p\nabla^{TM}_VV\rangle+\frac{1}{2}\langle
p\nabla^{TM}_VV,2\nabla^{TM}_UU\rangle-
\frac{1}{2\varepsilon}\langle
p\nabla^{TM}_VV,2\nabla^{TM}_UU\rangle+\langle
\nabla_V^{F^{\bot}}\nabla^{F^{\bot}}_U,V\rangle$$
$$=\frac{1}{2\varepsilon}\omega(p\nabla^{TM}_VV)(U,U)+O(1). \eqno(2.7)$$
Similarly, we have,\\
 $$\langle p\nabla^{TM,\varepsilon}_VU,\nabla^{TM}_UV+\nabla^{TM}_VU\rangle=\langle
 \nabla^{TM,\varepsilon}_VU,p(\nabla^{TM}_UV+\nabla^{TM}_VU)\rangle$$
 $$=\frac{1}{2\varepsilon}\langle
 p(\nabla^{TM}_UV+\nabla^{TM}_VU),\nabla^{TM}_UV+\nabla^{TM}_VU\rangle+O(1)$$
 $$=-\frac{1}{2\varepsilon}\omega(p(\nabla^{TM}_UV+\nabla^{TM}_VU))(U,V)+O(1).\eqno(2.8)$$
By (2.7), (2.8) and (2.34) in [LZ], we get,
$${\rm lim}_{\varepsilon\rightarrow
 0}\varepsilon\langle R^{TM,\varepsilon}(U,V)U,V\rangle=\frac{1}{4}\omega(p(\nabla^{TM}_UV+\nabla^{TM}_VU))(U,V)
-\frac{1}{2}\omega(p\nabla^{TM}_VV)(U,U).\eqno(2.9)$$
   By definition, $$-k^{TM,\varepsilon}=\sum_{i,j=1}^p\langle R^{TM,\varepsilon}(f_i,f_j)f_i,f_j\rangle
+\sum_{s,t=1}^q \varepsilon \langle
R^{TM,\varepsilon}(h_s,h_t)h_s,h_t\rangle$$
$$+2\sum_{i=1}^p\sum_{s=1}^q\langle R^{TM,\varepsilon}(f_i,h_s)f_i,h_s\rangle.\eqno(2.10)$$
  By (2.4), (2.6), (2.9), (2.10), then we get
 $${\rm lim}_{\varepsilon\rightarrow
0}(k^{TM,\varepsilon})=k^F+\Phi(\omega),\eqno(2.11)$$ where
$$\Phi(\omega)=\sum_{s,t=1}^q[-\frac{1}{4}\omega(p(\nabla^{TM}_{h_s}h_t+\nabla^{TM}_{h_t}h_s))(h_s,h_t)
+\frac{1}{2}\omega(p\nabla^{TM}_{h_t}h_t)(h_s,h_s)]$$
$$+\sum_{i=1}^q\sum_{s=1}^q\{
\frac{1}{2}\omega(p\nabla^{TM}_{f_i}f_i)(h_s,h_s)+\frac{1}{2}\omega(f_i)(f_i,p^{\bot}[f_i,h_s])$$
$$-\langle[f_i,A],h_s\rangle-\frac{1}{2}\omega(f_i)(h_s,A)\}.\eqno(2.12)$$
\indent Borrowing the idea in [LZ], we will prove a vanishing
theorem. By Theorem 1.1 in [LZ], we can write for $Y\in \Gamma
(F),U\in \Gamma (F^{\bot}),$
$$\nabla_Y^{F^{\bot},\varepsilon}U=
\widehat{\nabla_Y}U+O(\varepsilon)=\dot{\nabla_Y}U+\frac{\omega(Y)U}{2}+O(\varepsilon).\eqno(2.13)$$
then because the curvature of the Bott connection vanishes along
leaves, we have for $X,Y\in \Gamma (F),U,V\in \Gamma (F^{\bot}),$
$$\langle
R^{F^{\bot},\varepsilon}(X,Y)U,V\rangle=\langle{\widehat
R}^{F^{\bot}}(X,Y)U,V\rangle+O(\varepsilon)$$
$$=\langle
\dot{R}(X,Y)U,V\rangle+\frac{1}{2}\langle\{\dot{\nabla_X}\omega(Y)
-\dot{\nabla}_Y\omega(X)\}U,V\rangle$$ $$+
\frac{1}{4}\langle\{[\omega(X),\omega(Y)]-2\omega([X,Y])\}U,V\rangle+O(\varepsilon)$$
$$=\frac{1}{2}\langle\{\dot{\nabla_X}\omega(Y)-\dot{\nabla_Y}\omega(X)\}U,V\rangle+
\frac{1}{4}\langle\{[\omega(X),\omega(Y)]-2\omega([X,Y])\}U,V\rangle+O(\varepsilon).\eqno(2.14)$$
By $\nabla^{F^{\bot},\varepsilon}_VU=\nabla^{F^{\bot}}_VU$ and
(2.13), similar to (2.14), for $X\in \Gamma (F),U,V,Z,Z_1,Z_2\in
\Gamma (F^{\bot})$ we have
$$\langle
R^{F^{\bot},\varepsilon}(X,U)V,Z\rangle=O(1);\langle
R^{F^{\bot},\varepsilon} (U,V)Z_1,Z_2\rangle=O(1).\eqno(2.15)$$ So
by (2.14), (2.15), the sum of the last three terms in (2.30) in [LZ]
is
$$=\frac{1}{8}\sum_{i,j=1}^p\sum_{s,t=1}^q\langle{\widehat
R}^{F^{\bot}}(f_i,f_j)h_t,h_s\rangle
c(f_i)c(f_j)\widehat{c}(\sqrt{\varepsilon}h_s)\widehat{c}(\sqrt{\varepsilon}h_t)+O(\sqrt{\varepsilon}).\eqno(2.16)$$
In order to prove a vanishing theorem, we will estimate the norm of
$\widehat{c}(\sqrt{\varepsilon}h_t)$ under the metric
$g^{\varepsilon}$. By
$$||\widehat{c}(\sqrt{\varepsilon}h_t)||_{g^{\varepsilon}}=||\widehat{c}(h_t)||_g,\eqno(2.17)$$
so we do not consider the order of $\varepsilon$ from
$\widehat{c}(\sqrt{\varepsilon}h_t).$
 For other three terms including
$R^{\phi (F^{\bot}),\varepsilon}$ in (2.30) in [LZ], similar to
(2.16), their sum is
$$\frac{1}{2}\sum_{i,j=1}^pc(f_i)c(f_j)\widehat{R}^{\phi (F^{\bot})}(f_i,f_j)+O(\sqrt{\varepsilon}).\eqno(2.18)$$
Then by (2.30) in [LZ], (2.11),(2.16) and (2.18), we get
$$D^2_{F,\phi(F^{\bot}),\varepsilon}=-\triangle^{F,\phi(F^{\bot}),\varepsilon}+k^F
+\Phi(\omega)$$
$$+\frac{1}{8}\sum_{i,j=1}^p\sum_{s,t=1}^q\langle{\widehat
R}^{F^{\bot}}(f_i,f_j)h_t,h_s\rangle
c(f_i)c(f_j)\widehat{c}(\sqrt{\varepsilon}h_s)\widehat{c}(\sqrt{\varepsilon}h_t)$$
$$+\frac{1}{2}\sum_{i,j=1}^pc(f_i)c(f_j)\widehat{R}^{\phi
(F^{\bot})}(f_i,f_j)+O({\sqrt\varepsilon}).\eqno(2.19)$$ Let $x$ be
a point in $M$ and $||L||_x$ be the norm of the operator $L$ on
$(S(F)\widehat{\otimes}\wedge(F^{\bot,\star})\otimes\phi(F^{\bot}),g_x).$
 Let $$ A(\omega,\phi
(F^{\bot}))(x)=\frac{(k^F +\Phi(\omega))(x)}{4}
-||\frac{1}{2}\sum_{i,j=1}^pc(f_i)c(f_j)\widehat{R}^{\phi
(F^{\bot})}(f_i,f_j)||_x$$ $$
-||\frac{1}{8}\sum_{i,j=1}^p\sum_{s,t=1}^q{\widehat
R}^{F^{\bot}}\langle(f_i,f_j)h_t,h_s\rangle
c(f_i)c(f_j)\widehat{c}(h_s)\widehat{c}(h_t)||_x.\eqno(2.20)$$ By
(2.17), (2.19), (2.20), if $A(\omega,\phi (F^{\bot}))>0$, for
sufficiently small $\varepsilon>0$, we get
$$D^2_{F,\phi(F^{\bot}),\varepsilon}>0~ {\rm and}~\widehat{A}(F)L(F^{\bot}){\rm
ch}(\phi(F^{\bot}))=0.\eqno(2.21)$$ \noindent Let
$$\langle\widehat{A}(TM),[M]\rangle=\sum c_i\langle\widehat{A}(F)L(F^{\bot}){\rm
ch}(\phi_i(F^{\bot})),[M]\rangle,\eqno(2.22)$$ where $c_i$ is a
constant. By (2.21) and (2.22), we have\\

\noindent{\bf Theorem 2.1}~~{\it Let $(M,F,g^F)$ be a compact and
transversally oriented foliation with spin leave, if
$A(\omega,\phi_i (F^{\bot}))>0$ for any $\phi_i(F^{\bot})$ appeared
in (2.22), then
$<\widehat{A}(TM),[M]>=0.$}\\

\noindent {\bf Remark.} When $(M,F)$ is a Riemannian foliation, then
$A=\frac{k^F}{4}$, so we get the vanishing theorem of Connes. As in
[LZ], we can get the vanishing theorem of Connes for almost
Riemannian
foliations from Theorem 2.1.\\

\section{ A vanishing theorem when $F$ is not integrable}

\quad Recall that for a contact manifold $M$ (for definition, see
[Bl]), we have $TM=F\oplus F^{\bot}$ where $F^{\bot}$ is a line
bundle and $F$ is not integrable. In order to get a vanishing
theorem as in Section 2 for contact manifolds, we need to compute
the adiabatic limit of the scalar curvature in this case.\\
 \indent
Let $(M,g^{TM})$ be an oriented Riemannian manifold. Assume that
$$TM=F\oplus F^{\bot},~g^{TM}=g^F\oplus g^{F^{\bot}},\eqno(3.1)$$
where $F$ may not be integrable. Now we compute the adiabatic limit
of the scalar curvature in this case. We use the same notations as
in Section 2. Then we have similar formulas to (1.5)-(1.8) in [LZ].
For $\widehat{X}\in TM, X,Y,Z\in \Gamma(F),U,V\in \Gamma
(F^{\bot}),$
$$\langle\nabla^{F,\varepsilon}_{\widehat{X}}Y,Z\rangle_{g^{TM}}
=\langle\nabla^{F}_{\widehat{X}}Y,Z\rangle_{g^{TM}}+\frac{1}{2}(1-\frac{1}{\varepsilon})\langle
p^{\bot}[Y,Z],\widehat{X}\rangle_{g^{TM}}.\eqno(3.2)$$
 Especially, when $\widehat{X}=X\in \Gamma(F)$,we have
$$\nabla^{F,\varepsilon}_X=\nabla^F_X.\eqno(3.3)$$ Furthermore,
 $$\langle\nabla^{TM,\varepsilon}_{X}U,Y\rangle_{g^{TM}}
=\langle\nabla^{TM}_{X}U,Y\rangle_{g^{TM}}+\frac{1}{2}(1-\frac{1}{\varepsilon})\langle
[X,Y],U\rangle_{g^{TM}},\eqno(3.4)$$

$$\langle\nabla^{TM,\varepsilon}_{V}U,X\rangle_{g^{TM}}
=\langle\nabla^{TM}_{V}U,X\rangle_{g^{TM}}-\frac{1}{2}(1-\frac{1}{\varepsilon})\langle
X,\nabla^{TM}_VU+\nabla^{TM}_UV\rangle,\eqno(3.5)$$  and
$$\langle\nabla^{TM,\varepsilon}_{X}Y,U\rangle_{g^{TM}}
=\varepsilon\langle\nabla^{TM}_{X}Y,U\rangle_{g^{TM}}
+\frac{1}{2}(1-\varepsilon)\langle[X,Y],U\rangle_{g^{TM}}.
\eqno(3.6)$$ Especially when $X=Y$, we have
$$p^{\bot}\nabla^{TM,\varepsilon}_{X}X=\varepsilon
p^{\bot}\nabla^{TM}_{X}X.\eqno(3.7)$$ Furthermore,
$$\langle\nabla^{TM,\varepsilon}_{V}Y,U\rangle_{g^{TM}}
=-\frac{1}{2}\langle
Y,\nabla^{TM}_VU+\nabla^{TM}_UV\rangle+\frac{\varepsilon}{2}\langle
Y,[U,V]\rangle,\eqno(3.8)$$
$$\nabla^{F^{\bot},\varepsilon}_V=\nabla^{F^{\bot}}_V,\eqno(3.9)$$
 and
$$\langle\nabla^{F^{\bot},\varepsilon}_{X}U,V\rangle_{g^{TM}}
=\langle[X,U],V\rangle-\frac{1}{2}\langle
X,\nabla^{TM}_VU+\nabla^{TM}_UV\rangle-\frac{\varepsilon}{2}\langle
X,[U,V]\rangle.\eqno(3.10)$$ \indent As in Section 1 in [LZ], we can
still define the Bott connection $\dot{\nabla}$ which may not be
flat along leave, its dual connection $\dot{\nabla}^{\star}$ and
$\omega,\widehat{\nabla}.$ Then we still have for $X\in \Gamma(F)$
$${\rm lim}_{\varepsilon\rightarrow
0}\nabla^{F^{\bot},\varepsilon}_X=\widehat{\nabla}_X.\eqno(3.11)$$
In our case, we denote the scalar curvature (resp. curvature) by
$\overline{k}^{TM,\varepsilon}$ (resp.
$\overline{R}^{TM,\varepsilon}$) associated to $g^{\varepsilon}$. We
still denote the curvature by $R^{TM,\varepsilon}$ when $F$ is
integrable. That is, we just use the expression of
$R^{TM,\varepsilon}$ in [LZ]. By (2.32) in [LZ] and (3.2), (3.3),
(3.6), we have
$$\langle\overline{R}^{TM,\varepsilon}(f_i,f_j)f_i,f_j\rangle=
\langle R^{TM,\varepsilon}(f_i,f_j)f_i,f_j\rangle
+(-\frac{1}{2}-\frac{\varepsilon}{4}+\frac{3}{4\varepsilon})\langle
p^{\bot}[f_i,f_j],[f_i,f_j]\rangle_{g^{TM}}.\eqno(3.13)$$ By (2.33)
in [LZ] and (3.2), (3.7), (3.8), (3.9), (3.10), for $X\in \Gamma(F),
U\in\Gamma(F^{\bot}),$ we have
$$\langle\overline{R}^{TM,\varepsilon}(X,U)X,U\rangle=\langle
R^{TM,\varepsilon}(X,U)X,U\rangle+\frac{1}{2}(1-\varepsilon)\langle[X,p\nabla^{TM}_XU],U\rangle$$
$$
+\frac{1}{2}(1-\varepsilon)\langle[X,p\nabla^{TM,\varepsilon}_UX],U\rangle
-\frac{1}{2}(1-\varepsilon)\langle[p[X,U],X],U\rangle.\eqno(3.14)$$
By (2.34) in [LZ] and (3.8), (3.9), (3.10), for
$U,V\in\Gamma(F^{\bot})$, we have
$$\langle\overline{R}^{TM,\varepsilon}(U,V)U,V\rangle=\langle
R^{TM,\varepsilon}(U,V)U,V\rangle.\eqno(3.15)$$ By (2.32) in [LZ]
and (3.13), when $\varepsilon\rightarrow 0$, we have
$$\langle\overline{R}^{TM,\varepsilon}(f_i,f_j)f_i,f_j\rangle\sim
\langle R^{F}(f_i,f_j)f_i,f_j\rangle-\frac{1}{2}\langle
p^{\bot}[f_i,f_j],[f_i,f_j]\rangle_{g^{TM}}+O(\frac{1}{\varepsilon})\eqno(3.16)$$
and
$${\rm lim}_{\varepsilon\rightarrow 0}\varepsilon\langle\overline{R}^{TM,\varepsilon}(f_i,f_j)f_i,f_j\rangle
=\frac{3}{4}\langle
p^{\bot}[f_i,f_j],[f_i,f_j]\rangle_{g^{TM}}.\eqno(3.17)$$ So
$$\sum_{1\leq i,j\leq p}{\rm lim}_{\varepsilon\rightarrow
0}\varepsilon\langle\overline{R}^{TM,\varepsilon}(f_i,f_j)f_i,f_j\rangle
=\frac{3}{4}\sum_{1\leq i,j\leq
p}||p^{\bot}[f_i,f_j]||^2,\eqno(3.18)$$ which is globally defined
and is positive when $F$ is not integrable. By (2.33), (2.35) in
[LZ] ((2.35) is still correct in this case) and (3.14), we have
$${\rm lim}_{\varepsilon\rightarrow 0}\varepsilon\langle\overline{R}^{TM,\varepsilon}(X,U)X,U\rangle
$$
$$={\rm lim}_{\varepsilon\rightarrow
0}\left[\varepsilon\langle
R^{TM,\varepsilon}(X,U)X,U\rangle+\frac{\varepsilon}{2}(1-\varepsilon)
\langle[X,p\nabla^{TM,\varepsilon}_UX],U\rangle\right]$$
$$={\rm lim}_{\varepsilon\rightarrow
0}\left[\frac{\varepsilon^2}{2}\langle
X,[U,p^{\bot}\nabla^{TM,\varepsilon}_UX]\rangle+\varepsilon\langle
[X,p^{\bot}\nabla_U^{TM,\varepsilon}X],U\rangle\right.$$
$$\left.-\frac{\varepsilon}{2}\langle
X,\nabla^{TM}_{p^{\bot}\nabla^{TM,\varepsilon}_UX}U+\nabla^{TM}_U{p^{\bot}\nabla^{TM,\varepsilon}_UX}\rangle
+\frac{\varepsilon}{2}(1-\varepsilon)
\langle[X,p\nabla^{TM,\varepsilon}_UX],U\rangle\right]$$
$$={\rm lim}_{\varepsilon\rightarrow
0}\frac{\varepsilon}{2}(1-\varepsilon)
\langle[X,p\nabla^{TM,\varepsilon}_UX],U\rangle\eqno(3.19)$$ By
(3.2), we get $$\langle[X,p\nabla^{TM,\varepsilon}_UX],U\rangle
=\langle\nabla^{TM}_UX,p\nabla^{TM}_XU\rangle+\frac{1}{2}(1-\frac{1}{\varepsilon})\langle
[X,p\nabla^{TM}_XU],U\rangle-\langle\nabla^{TM}_{p\nabla^{TM,\varepsilon}_UX}X,U\rangle\eqno(3.20)$$
\noindent and
$$p\nabla^{TM,\varepsilon}_UX=\sum_{j=1}^p\langle\nabla^{TM,\varepsilon}_UX,f_j\rangle
f_j=p\nabla^{TM}_UX+\frac{1}{2}(1-\frac{1}{\varepsilon})\sum_{j=1}^p\langle[X,f_j],U\rangle
f_j.\eqno(3.21)$$ \noindent
 by (3.19-3.21), we have
$${\rm lim}_{\varepsilon\rightarrow 0}\varepsilon\langle\overline{R}^{TM,\varepsilon}(X,U)X,U\rangle
=-\frac{1}{4}\langle[X,p\nabla^{TM}_XU],U\rangle+\frac{1}{4}\langle\nabla^{TM}_{\sum_{j=1}^p\langle[X,f_j],U\rangle
f_j}X,U\rangle.\eqno(3.22)$$ \noindent By (2.34) in [LZ] and (3.15),
we get
$${\rm lim}_{\varepsilon\rightarrow 0}\varepsilon^2\langle\overline{R}^{TM,\varepsilon}(U,V)U,V\rangle
=0\eqno(3.23)$$ By (3.17), (3.22) and (3.23), we have
$$4B:={\rm lim}_{\varepsilon\rightarrow
0}\varepsilon\overline{k}^{TM,\varepsilon}=- \frac{3}{4}\sum_{1\leq
i,j\leq p}||p^{\bot}[f_i,f_j]||^2$$
$$-\frac{1}{2}\sum_{i=1}^p\sum_{s=1}^q\{
\nabla^{TM}_{\sum_{j=1}^p\langle[f_i,f_j],h_s\rangle
f_j}f_i,h_s\rangle-\langle[f_i,p\nabla^{TM}_{f_i}h_s],h_s\rangle\}$$
$$=- \frac{3}{4}\sum_{1\leq
i,j\leq
p}||p^{\bot}[f_i,f_j]||^2-\frac{1}{2}\sum_{i=1}^p\sum_{s=1}^q
||p\nabla^{TM}_{f_i}h_s||^2+\frac{1}{2}\sum_{i=1}^p\sum_{j=1}^p
||p^{\bot}\nabla^{TM}_{f_j}f_i||^2,\eqno(3.24)$$ \noindent which is
globally defined and vanishes when $F$ is integrable. By (2.30) in
[LZ] and (3.11), similar to Section 2, we have
$$\varepsilon D^2_{F,\phi(F^{\bot}),\varepsilon}=-\varepsilon\Delta^{F,\phi(F^{\bot}),\varepsilon}
+\frac{\varepsilon\overline{k}^{TM,\varepsilon}}{4}+O(\varepsilon).\eqno(3.25)$$
By (3.24) and (3.25) and that
$-\varepsilon\Delta^{F,\phi(F^{\bot}),\varepsilon}$ is nonnegative,
we get for sufficient small $\varepsilon>0,$ if $B>0,$ then
$D^2_{F,\phi(F^{\bot}),\varepsilon}>0.$ So similarly to Theorem 2.1,
we have\\

 \noindent{\bf Theorem 3.1}~~{\it Let $M$ be a compact
oriented manifold, if 1) $TM=F\oplus F^{\bot}$ and $F$ is not
integrable and oriented spin, 2) $B>0$, then
$<\widehat{A}(TM),[M]>=0.$}\\

\indent For many cases, $TM=F^{\bot}\oplus F$ and $F^{\bot}$ is not
integrable. When $M$ be an compact contact metric manifold and
dim$M\geq 3$, then $M$ has a canonical ${\rm spin^{c}}$ structure
(for details, see [Pe]). In this case, $TM=L\oplus L^{\bot}$ and
$L^{\bot}$ is not integrable. Let
$g^{TM,\varepsilon}=g^L\oplus\frac{1}{\varepsilon}g^{L^{\bot}}.$
 Let $D^{\widetilde{L}}$ be the ${\rm spin^{c}}$ Dirac operator associated to the
 Levi-Civita connection and the general complex determine line bundle
 $\widetilde{L}.$ By [LM], we have
 $$(D^{\widetilde{L},\varepsilon})^2=-\Delta^{\widetilde{L},\varepsilon}+\frac{\overline{k}^{TM,\varepsilon}}{4}
 +\frac{\sqrt{-1}}{2}\Omega^{\widetilde{L}}.\eqno(3.26)$$
 By (3.26), similar to the proof of Theorem 3.1, we have\\

 \noindent{\bf Theorem 3.2}~~{\it For the above
 $M$ and $~D^{\widetilde{L}}$, if
 $B>0,$ then ${\rm ker}D^{\widetilde{L},\varepsilon}=0$ for sufficient small
 $\varepsilon>0.\\$}

Similarly, let $M^{2n}$ be an compact oriented manifold. We assume
$TM=F\oplus F^{\bot}$ and $F$ is not integrable. Let
$D_{\varepsilon}=d_{\varepsilon}+\delta_{\varepsilon}$ be the
de-Rham Hodge operator (resp. signature operator) associated to
$(M,g^{\varepsilon}).$ If $B>0$, then the Euler number (resp.
signature) of $M$ is zero.\\

 \section{A vanishing theorem for complex foliations}

 \quad Let $M$ be an compact connected complex manifold of complex
 dimension $n$ with a complex foliated structure of complex dimension
 $p$. That is, $M$ is the disjoint union of its complex submanifolds
 of complex dimension $p$ which locally are defined by
 $dz_{p+1}=\cdots=dz_{p+q}=0,$ where $p+q=n$ and
 ${z_1=x_1+iy_1,\cdots,z_n=x_n+iy_n}$ is the complex coordinate of
 $M.$ We consider $M$ as an almost  complex manifold with the
 canonical almost complex structure $J:$
 $$J(\frac{\partial}{\partial x_k})=\frac{\partial}{\partial y_k};~
J(\frac{\partial}{\partial y_k})=-\frac{\partial}{\partial
x_k}.\eqno(4.1)$$ Let $TM$ be the holomorphic tangent bundle on $M$
and let $T_{\bf R}M$ be the real tangent bundle of $M$ as an real
manifold. Let $F$ be the real tangent bundle of complex leave and
locally
$$F={\rm span}_{\bf R}\{\frac{\partial}{\partial x_1},\frac{\partial}{\partial
y_1},\cdots,\frac{\partial}{\partial x_p},\frac{\partial}{\partial
y_p}\},\eqno(4.2)$$ then $J|_F:F\rightarrow F$ is a complex
structure of the bundle $F$. As in Section 2, we take $g^{T_{\bf
R}M}=g^F\oplus g^{F^{\bot}}$ and $T_{\bf R}M/F\cong F^{\bot}.$ Since
$J:F\oplus F^{\bot}\rightarrow F\oplus F^{\bot}$ and
$J|_F:F\rightarrow F$ are isomorphism, then $J:F^{\bot}\rightarrow
F^{\bot}$ is an isomorphism, so $J$ is also a complex structure of
$F^{\bot}.$ We take a positive definite Hermitian structure
$(F,J,H^{F})$ (resp. $(F^{\bot},J,H^{F^{\bot}})$) of $(F,J)$ (resp.
$(F^{\bot},J)$) such that $$g^F={\rm Re}H^F;~ g^{F^{\bot}}={\rm
Re}H^{F^{\bot}}.\eqno(4.3)$$ Here ${\rm Re}H$ (resp. ${\rm Im}H$)
denotes the real part (resp. imaginary part) of the Hermitian metric
$H$. Let $\omega_1$ (resp. $\omega_2$) be the K\"{a}hler form
associated to $H^F$ (resp. $H^{F^{\bot}}$). Then for $X_1,Y_1\in
F,~X_2,Y_2\in F^{\bot},$
$$\omega_1(X_1,Y_1)=g^F(X_1,JY_1);~\omega_2(X_2,Y_2)=g^{F^{\bot}}(X_2,JY_2).
\eqno(4.4)$$ \noindent For any $\varepsilon>0$, we define a positive
definite Hermitian structure $H^{\varepsilon}$ on $T_{\bf R}M$ by
$$H^{\varepsilon}=H^F\oplus\frac{1}{\varepsilon}H^{F^{\bot}},\eqno(4.5)$$
 Then
$$g^{T_{\bf R}M,\varepsilon}={\rm
Re}H^{\varepsilon};~\widehat{\omega}^{\varepsilon}=\omega_1\oplus
\frac{1}{\varepsilon}\omega_2,\eqno(4.6)$$ \noindent and we write
$\widehat{\omega}=\widehat{\omega}^1.$\\
 \indent Now, we recall the
Bismut-Kodaira-Nakano formula (for details, see Section 2 in [Bi]).
Let $\xi$ be a holomorphic Hermitian vector bundle of complex
dimension $l.$ Let $\nabla^{\xi}$ be the holomorphic Hermitian
connection on $\xi$, whose curvature is denoted by
$(\nabla^{\xi})^2$. $\wedge T^{\star (0,1)}M$ denotes the algebra of
forms of type $(0,p)~ (0\leq p\leq n).$ Let $\overline{\partial}$ be
the Dolbeault operator acting on the set $\Gamma$ of smooth sections
of $\wedge T^{\star (0,1)}M\otimes \xi$ equipped with the natural
$L^2$-Hermitian product. Let $\overline{\partial}^{\star}$ be the
formal adjoint of $\overline{\partial}$ with respect to this
Hermitian metric. Let $\nabla^{TM}$ be the holomorphic Hermitian
connection on $TM$ associated to the hermitian metric $g^{TM}$
induced by $g^{T_{\bf R}M}=g^F\oplus g^{F^{\bot}}.$ Let $\omega$ be
the K\"{a}hler form associated to $H=H^F\oplus H^{F^{\bot}}.$
Locally,
$$\wedge T^{\star (0,1)}M\otimes \xi=S(T_{\bf
R}M)\otimes(\lambda\otimes\xi),\eqno(4.7)$$ \noindent where
$S(T_{\bf R}M)$ denotes spinors bundle and $\lambda$ denotes the
square root of ${\rm det}(TM).$ Then the Levi-Civita connection
$\nabla^L$ associated to $g^{T_{\bf R}M}$ has a lift on $S(T_{\bf
R}M)$, still denoted by $\nabla^L.$ Then by Theorem 2.3 in [Bi], we
have:
$$2(\overline{\partial}+\overline{\partial}^{\star})^2=-\sum_{i=1}^{2n}
((\nabla^E_{e_i})^2-\nabla^E_{\nabla^L_{e_i}e_i})^2+\frac{k}{4}$$
$$+^c((\nabla^{\xi})^2+\frac{1}{2}{\rm Tr}[(\nabla^{TM})^2]\otimes
I_{\xi})-\frac{\sqrt{-1}}{2}
{^c(\overline{\partial}\partial\widehat{\omega})}-
\frac{1}{8}||(\partial-\overline{\partial})\widehat{\omega}||^2,\eqno(4.8)$$
where
$E=-\frac{\sqrt{-1}}{4}(\partial-\overline{\partial})\widehat{\omega}$
and $\nabla^E=\nabla^L+S^E$ and $\langle S^E(X)Y,Z\rangle=2E(X,Y,Z)$
for $X, Y,Z\in \Gamma (T_{\bf R}M)$ and $\{e_1,\cdots,e_{2n}\}$ is
an orthonormal basis associated to $(T_{\bf R}M,g).$ If we replace
$g$ by $g^{\varepsilon},$ then by (4.8), we have
$$2(\overline{\partial}+\overline{\partial}^{\star}_{\varepsilon})^2
=-\Delta^E_{\varepsilon} +\frac{k^{\varepsilon}}{4}
+^c((\nabla^{\xi})^2+\frac{1}{2}{\rm
Tr}[(\nabla^{TM,\varepsilon})^2]\otimes I_{\xi})_{\varepsilon}$$ $$-
\frac{\sqrt{-1}}{2}
{^c(\overline{\partial}\partial\widehat{\omega}^{\varepsilon})}
-\frac{1}{8}||(\partial-\overline{\partial})\widehat{\omega}^{\varepsilon}||
^2_{g^{\varepsilon}}.\eqno(4.9)$$ Next we will compute the adiabatic
limits of some terms in (4.9). We assume
$$\omega_2=\sum_{1\leq i<j\leq 2n} f_{i,j}\overline{dx_i}\wedge \overline{dx_j}
;~f_{i,j}= \omega_2(\frac{\partial}
{\partial\overline{x_i}},\frac{\partial}{\partial\overline{x_j}}),\eqno(4.10)$$
where
$\{d\overline{x_1},\cdots,d\overline{x_{2n}}\}=\{dx_1,dy_1,\cdots,dx_n,dy_n\}.$
By (4.2) and $\omega_2(X,Y)=0$ if $X$ or $Y$ in $F$, we know
$$f_{i,j}=0,~{\rm when}~ i~ {\rm or}~ j\leq 2p.\eqno(4.11)$$
By direct computations and (4.11), we get
$$(\partial-\overline{\partial})\omega_2(f_i,f_j,f_k)=(\partial-\overline{\partial})\omega_2(f_i,f_j,h_k)=0;
\eqno(4.12)$$
$$\overline{\partial}\partial\omega_2(f_i,f_j,f_k,f_l)=\overline{\partial}\partial\omega_2(f_i,f_j,f_k,h_l)=0.
\eqno(4.13)$$ By (4.6) and (4.13), we have
$$^c(\overline{\partial}\partial\widehat{\omega}^{\varepsilon})=\sum_{1\leq
i<j<k<l\leq
2p}(\overline{\partial}\partial\omega_1)(f_i,f_j,f_k,f_l)c(f_i)c(f_j)c(f_k)c(f_l)$$
$$+ \sum_{1\leq i<j\leq 2p}\sum_{1\leq k<l\leq
2q}(\overline{\partial}\partial\omega_2)(f_i,f_j,h_k,h_l)c(f_i)c(f_j)c(\sqrt{\varepsilon}h_k)
c(\sqrt{\varepsilon}h_l)+O(\sqrt{\varepsilon}).\eqno(4.14)$$ By
(4.6) (4.12) and the definition of the norm, we get
$$||(\partial-\overline{\partial})\widehat{\omega}^{\varepsilon}||
^2_{g^{\varepsilon}}=\sum_{1\leq i<j<k\leq
2p}[\sqrt{-1}(\partial-\overline{\partial})\omega_1(f_i,f_j,f_k)]^2$$
$$+ \sum_{1\leq i\leq 2p}\sum_{1\leq j<k\leq
2q}[\sqrt{-1}(\partial-\overline{\partial})\omega_2(f_i,h_j,h_k)]^2+O(\varepsilon).\eqno(4.15)$$
\noindent Since $^c(\nabla^{\xi})^2$ is independent of
$g^{\varepsilon}$, we have
$$^c(\nabla^{\xi})^2=\sum_{1\leq i<j\leq
2p}(\nabla^{\xi})^2(f_i,f_j)c(f_i)c(f_j)+O(\sqrt{\varepsilon}).\eqno(4.16)$$
Next we consider $^c({\rm Tr}[(\nabla^{TM,\varepsilon})^2]).$ We
extend $g^{T_{\bf R}M}$ (resp. $g^{T_{\bf R}M,\varepsilon}$) to an
Hermitian metric $g^{T_{\bf R}M\otimes {\bf c}}$ (resp. $g^{T_{\bf
R}M\otimes {\bf C},\varepsilon}$), then we get an Hermitian metric
$g^{TM}$ (resp. $g^{TM,\varepsilon}$) by restricting the $g^{T_{\bf
R}M\otimes {\bf C}}$ (resp. $g^{T_{\bf R}M\otimes {\bf
C},\varepsilon}$) to $TM.$ Let $\{ \frac{\partial}{\partial
z_1},\cdots,\frac{\partial}{\partial z_n}\}$ be a holomorphic local
basis of $TM$ and
$$\widehat{H}=(H_{\alpha\beta})_{n\times
n};~\widehat{H}_{\alpha\beta}=g^{TM}(\frac{\partial}{\partial
z_{\alpha}},\frac{\partial}{\partial z_{\beta}}).\eqno(4.17)$$
\noindent Let $\widetilde{\omega}$ (resp. $\Omega$) be the
connection (resp. curvature) matrix associated to the Hermitian
holomorphic connection $(TM,\nabla^{TM},g^{TM})$ under the basis $\{
\frac{\partial}{\partial z_1},\cdots,\frac{\partial}{\partial
z_n}\}.$ Then
$$\widetilde{\omega}=\partial\widehat{H}.\widehat{H}^{-1};~\Omega=d\widetilde{\omega}
-\widetilde{\omega}\wedge\widetilde{\omega}.\eqno(4.18)$$ By $T_{\bf
R}M\otimes{\bf C}=F\otimes{\bf C}\oplus F^{\bot}\otimes{\bf C}$ and
$(M,F)$ is a complex foliation, we know
$$TM=T^{1,0}F\oplus T^{1,0}F^{\bot},\eqno(4.19)$$
where locally
$$T^{1,0}F={\rm span}_{\bf C}\{\frac{\partial}{\partial z_1},\cdots,\frac{\partial}{\partial
z_p}\}\eqno(4.20)$$ and $T^{1,0}F^{\bot}$ is the orthogonal
complementary bundle of $T^{1,0}F$ in $TM.$ We let
$$\frac{\partial}{\partial z_r}=\sum_{j=1}^pa_{rj}\frac{\partial}{\partial z_j}
+\sum_{s=1}^qb_{rs}\overline{e_s},~p+1\leq r\leq n,\eqno(4.21)$$
\noindent where $\{\overline{e_1},\cdots,\overline{e_q}\}$ is a
complex orthonormal basis of $T^{1,0}F^{\bot}$. Let
$B=(b_{rs})_{q\times q}$, then $$ \left[\begin{array}{lcr}
  \ \frac{\partial}{\partial z_1} \\
    \  \vdots \\
    \ \frac{\partial}{\partial z_n}
\end{array}\right]
=\left[\begin{array}{lcr}
  \ I_{p\times p} &  0 \\
    \ \star  & B
\end{array}\right]
\left[\begin{array}{lcr}
  \ \frac{\partial}{\partial z_1} \\
    \  \vdots \\
    \ \frac{\partial}{\partial z_p} \\
    \ \overline{e_1} \\
    \ \vdots \\
    \ \overline{e_s}
\end{array}\right],\eqno(4.22)$$
\noindent where $B$ is invertible. Let
$$\widehat{H}=\left[\begin{array}{lcr}
  \ H_{pp} &  H_{pq} \\
    \ H_{qp}  & H_{qq}
\end{array}\right],
\eqno(4.23)$$ \noindent where
$H_{pp}=[\widehat{H}_{\alpha\beta}]_{1\leq \alpha,\beta\leq p}$ is
invertible. By (4.20),(4.21) and
$$g^{T_{\bf R}M\otimes {\bf c},\varepsilon}=g^{F\otimes {\bf
C}}\oplus \frac{1}{\varepsilon}g^{F^{\bot}\otimes {\bf
C}},\eqno(4.24)$$ \noindent then under the metric
$g^{TM,\varepsilon},$ we have
$$\widehat{H}^{\varepsilon}=\left[\begin{array}{lcr}
  \ H_{pp} &  H_{pq} \\
    \ H_{qp}  &
    \overline{H_{qq}}+\frac{1}{\varepsilon}\widehat{H_{qq}}
\end{array}\right],
\eqno(4.25)$$ \noindent where $\widehat{H_{qq}}=B\overline{B}^t$
which is invertible. By the definition of the inverse of $
\widehat{H}^{\varepsilon},$ we get
$$(\widehat{H}^{\varepsilon})^{-1}=\frac{1}{{\rm det}H_{pp}{\rm
det}\widehat{H}_{qq}+O(\varepsilon)} \left[\begin{array}{lcr}
  \ A_{11}+O(\varepsilon) &  O(\varepsilon) \\
    \ A_{21}\varepsilon+O(\varepsilon^2)  & A_{22}\varepsilon+O(\varepsilon^2)
   \end{array}\right],
\eqno(4.26)$$ \noindent where
$$A_{11}={\rm det}H_{pp}{\rm
det}\widehat{H}_{qq}H_{pp}^{-1};~A_{22}={\rm det}H_{pp}{\rm
det}\widehat{H}_{qq}\widehat{H_{qq}}^{-1}.\eqno(4.27)$$ \noindent By
(4.18),(4.25) and (4.26),
$$\widetilde{\omega}^{\varepsilon}=\frac{1}{{\rm det}H_{pp}{\rm
det}\widehat{H}_{qq}} \left[\begin{array}{lcr}
  \ \partial H_{pp}.A_{11} &  0 \\
    \ \partial H_{qp}.A_{11}+\partial\widehat{H_{qq}}.A_{21}  & \partial\widehat{H_{qq}}.A_{22}
   \end{array}\right]
+O(\varepsilon).\eqno(4.28)$$ \noindent Then by
(4.18),(4.27),(4.28), we get
$$
{\rm Tr}[(\nabla^{TM,\varepsilon})^2]={\rm
Tr}[\Omega^{\varepsilon}]={\rm Tr}[\Omega^{T^{1,0}F}]+{\rm
Tr}[\Omega^{T^{1,0}F^{\bot}}]+O(\varepsilon),\eqno(4.29)$$ \noindent
where
$$\Omega^{T^{1,0}F}=d(\partial H_{pp}.H_{pp}^{-1})-(\partial
H_{pp}.H_{pp}^{-1})\wedge (\partial H_{pp}.H_{pp}^{-1});$$
$$\Omega^{T^{1,0}F^{\bot}}=d(\partial \widehat{H}_{qq}.\widehat{H}_{qq}^{-1})-(\partial
\widehat{H}_{qq}.\widehat{H}_{qq}^{-1})\wedge (\partial
\widehat{H}_{qq}.\widehat{H}_{qq}^{-1}).\eqno(4.30)$$ \noindent By
(4.29), we obtain $$^c\left({\rm
Tr}[(\nabla^{TM,\varepsilon})^2]\right) =\sum_{1\leq i<j\leq
2p}\left\{{\rm Tr}[\Omega^{T^{1,0}F}(f_i,f_j)]+ {\rm
Tr}[\Omega^{T^{1,0}F^{\bot}}(f_i,f_j)]\right\}c(f_i)c(f_j)+O(\sqrt{\varepsilon}).\eqno(4.31)$$
\noindent {\bf Remark.} In fact the term $O(\sqrt{\varepsilon})$ in
(4.29) and (4.31) vanishes. We let
$${B_0}=\left[\begin{array}{lcr}
  \ I_{p\times p} &  0 \\
    \ \star  & B
\end{array}\right],\eqno (4.32)$$
\noindent then
$$\widehat{H}^{\varepsilon}={B_0}\left[\begin{array}{lcr}
  \ H_{pp} &  0 \\
    \ 0  & \frac{1}{\varepsilon}
\end{array}\right]\overline{B_0}^t;~(\widehat{H}^{\varepsilon})^{-1}=\overline{B_0^{-1}}^t
\left[\begin{array}{lcr}
  \ H_{pp}^{-1} &  0 \\
    \ 0  & \varepsilon
\end{array}\right]{{B_0}}^{-1}.\eqno(4.33)$$
\noindent By (4.18), (4.33) and the trace property, we have
$${\rm Tr}(\widetilde{\omega}^{\varepsilon})={\rm Tr}[\partial
H_{pp}.H_{pp}^{-1}+\partial {B_0}{B_0}^{-1}+\partial
\overline{B_0}\overline{B_0}^{-1}].\eqno(4.34)$$ \noindent So the
vanishing of the $O(\sqrt{\varepsilon})$ in (4.29) and (4.31) comes
from (4.34) and
$${\rm
Tr}[\Omega^{\varepsilon}]=\overline{\partial}{\rm
Tr}(\widetilde{\omega}^{\varepsilon}).\eqno(4.35)$$\\
\indent Let $\phi_{ij}$ be the $(i,j)$-component of $\phi\in \wedge
^pT^{\bf{R},\star}M\otimes {\bf C}$ associated to the decomposition
$$\wedge ^pT^{\bf{R},\star}M\otimes {\bf
C}=\oplus_{i+j=p}\wedge^iF^{\star}\otimes\wedge^j{F^{\bot,\star}}\otimes{\bf{C}}.\eqno(4.36)$$
\noindent We set
$$\Psi(x)=\frac{(k^F +\Phi(\omega))(x)}{4}-||\frac{\sqrt{-1}}{2}
 {^c[(\overline{\partial}\partial\omega_2)_{2,2}]}||_x$$
$$-||\frac{\sqrt{-1}}{2}
{^c[(\overline{\partial}\partial\omega_1)_{4,0}]}||_x
-\frac{1}{8}||\sqrt{-1}[(\partial-\overline{\partial})\omega_2]_{1,2}||^2_x$$
$$-\frac{1}{8}||\sqrt{-1}[(\partial-\overline{\partial})\omega_1]_{3,0}||^2_x
-||^c[(\nabla^{\xi})^2]_{2,0}||_x-\frac{1}{2}||^c[{\rm
Tr}(\Omega^{T^{1,0}F}+\Omega^{T^{1,0}F^{\bot}})]_{2,0}||_x,\eqno(4.37)$$
\noindent where $||.||_x$ denotes the norm of a linear operator
acting on $\left((\wedge T^{\star (0,1)}M\otimes \xi)_x,g_x\right).$
As in the
discussions in Section 2, we have\\

 \noindent{\bf Theorem 4.1}~~{\it Let $(M,F)$ be a compact complex
 foliation. If $\Psi>0,$ then the Euler number ${\rm Eul}(\xi)$ of the holomorphic
 bundle $\xi$ vanishes.}\\

By Theorem 2.11 in [Bi], we have:\\

 \noindent{\bf Corollary 4.2}~~{\it Let $(M,F)$ be a compact complex
 foliation. If $\Psi>0$ and
 $\overline{\partial}\partial(\omega_1+\omega_2)=0$, then there
 exists a $(2n-1)$-form $\tau$ such that}
 $$\left\{\widehat{A}\left(\frac{R^{-E}}{2\pi}\right){\rm
 exp}\left(-\frac{1}{2\sqrt{-1}\pi}{\rm
 Tr}\left[\frac{(\nabla^{TM})^2}{2}\right]\right){\rm Tr}\left[{\rm
 exp}\left(-\frac{(\nabla^{\xi})^2}{2\sqrt{-1}\pi}\right)\right]\right\}^{\rm
 max}=d\tau.\eqno(4.38)$$

\section{A Kastler-Kalau-Walze type theorem for
   foliations }

\quad Several years ago, Connes made a challanging observation that
the noncommutative residue (see [FGV] or [Wo]) of the $-2l+2$ power
of the Dirac operator on $2l~(l\geq 2)$-dimensional spin manifolds
was proportional to the Einstein-Hilbert action, which was called
Kastler-Kalau-Walze Theorem now. In [Ka], Kastler gave a brute-force
proof of this theorem. In [KW], Kalau and Walze proved this theorem
by the normal coordinates way simultaneously.
In [Ac], Ackermann gave a note on a new proof of this theorem by using heat kernel expansion.\\
  \indent Let $(M^n,F)$ be a compact even dimensional foliation. In
  [Ko1], the author defined the tangential pseudodifferential
   operator algebra for foliations and a noncommutative residue on
   it. We assume $M$ has spin leave. In order to give a Kastler-Kalau-Walze type theorem for
   foliations, it is natural to consider the noncommutative residue
   in [Ko1]
   of the $-n+2$
power of the Dirac operator along leave which is in the tangential
pseudodifferential
   operator algebra. But the computations seems to be a little
   complicated , which comes from the tangential
pseudodifferential calculas. Borrowing the idea in [LZ], we consider
the sub-Dirac operator and the noncommutative residue on the
classical pseudodifferential
   operator algebra on $M.$ That is, we will compute
$${\rm lim}_{\varepsilon\rightarrow 0}{\rm
Res}(D^{-n+2}_{F,\varepsilon}),$$ \noindent
where $D_{F,\varepsilon}=D_{F,\phi(F^{\bot})\varepsilon}|_{\phi=1}.$\\
\indent By (2.20) in [LZ], then
$$D^2_F=-\triangle^{F}+\frac{k^{TM}}{4}+Q\eqno(5.1)$$
\noindent where
$$Q=\frac{1}{4}\sum_{i=1}^p\sum_{r,s,t=1}^q\langle
R^{F^{\bot}}(f_i,h_r)h_t,h_s\rangle
c(f_i)c(h_r)\widehat{c}(h_s)\widehat{c}(h_t)$$
$$+\frac{1}{8}\sum_{i,j=1}^p\sum_{s,t=1}^q\langle
R^{F^{\bot}}(f_i,f_j)h_t,h_s\rangle
c(f_i)c(f_j)\widehat{c}(h_s)\widehat{c}(h_t)$$
$$
+\frac{1}{8}\sum_{s,t,r,l=1}^q\langle
R^{F^{\bot}}(h_r,h_l)h_t,h_s\rangle
c(h_r)c(h_l)\widehat{c}(h_s)\widehat{c}(h_t). \eqno(5.2)$$ \noindent
By [Ac], we know that:
$${\rm
Res}(D^{-n+2}_{F})=c_0\int_M{\rm
Tr}_{S(F)\widehat{\otimes}\wedge(F^{\bot,\star})}(-\frac{k^{TM}}{12}-Q)d{\rm
vol}_g,\eqno(5.3)$$ where $c_0=\frac{2}{(\frac{n}{2}-2)!\times
(4\pi)^{\frac{n}{2}}}.$ By the identity
$$d{\rm
vol}_{g^{\varepsilon}}=\frac{1}{\varepsilon^{\frac{q}{2}}}d{\rm
vol}_g,\eqno(5.4)$$ \noindent we get
$${\rm lim}_{\varepsilon\rightarrow 0}\varepsilon^{\frac{q}{2}}{\rm
Res}(D^{-n+2}_{F,\varepsilon}) =c_0\int_M{\rm
lim}_{\varepsilon\rightarrow 0} {\rm
Tr}_{(S(F)\widehat{\otimes}\wedge(F^{\bot,\star}),g^{\varepsilon})}\left(
-\frac{k^{TM,\varepsilon}}{12}-Q^{\varepsilon}\right)\varepsilon^{\frac{q}{2}}
d{\rm vol}_{g^{\varepsilon}}$$
$$=c_0\int_M{\rm
lim}_{\varepsilon\rightarrow 0} {\rm Tr}_{g^{\varepsilon}}\left(
-\frac{k^{TM,\varepsilon}}{12}-Q^{\varepsilon}\right) d{\rm
vol}_{g}.\eqno(5.5)$$ By (2.14) and ${\rm
Tr}_{g^{\varepsilon}}(\widehat{c}(\sqrt{\varepsilon}h_s)) ={\rm
Tr}_{g}(\widehat{c}(h_s)),$ we have $${\rm
lim}_{\varepsilon\rightarrow 0}{\rm
Tr}_{g^{\varepsilon}}Q^{\varepsilon} ={\rm Tr} [
\frac{1}{8}\sum_{i,j=1}^p\sum_{s,t=1}^q\langle
\widehat{R}^{F^{\bot}}(f_i,f_j)h_t,h_s\rangle
c(f_i)c(f_j)\widehat{c}(h_s)\widehat{c}(h_t)]=0.\eqno(5.6)$$
\noindent In the last equality of (5.6), we have used the identity
$${\rm Tr}[c(f_i)c(f_j)]={\rm
Tr}[\widehat{c}(h_s)\widehat{c}(h_t)]=0,~{\rm for }~ i\neq j,~s\neq
t.\eqno(5.7)$$ \noindent By (2.11) (5.5) (5.6), we have
$${\rm lim}_{\varepsilon\rightarrow 0}\varepsilon^{\frac{q}{2}}{\rm
Res}(D^{-n+2}_{F,\varepsilon})
=\widehat{c_0}\int_M[k^F+\Phi(\omega)]d{\rm vol}_{g},\eqno(5.8)$$
\noindent where $\widehat{c_0}=-\frac{c_0}{12}{\rm
Rk}[S(F)\widehat{\otimes}\wedge(F^{\bot,\star})]$ and ${\rm
Rk}[S(F)\widehat{\otimes}\wedge(F^{\bot,\star})]$ equals
$2^{\frac{p}{2}+q}$ (resp. $2^{\frac{p-1}{2}+q}$) when $p$ is even
(resp. odd). So we get\\

\noindent{\bf Theorem 5.1}~~{\it Let $(M^n,F)$ be a compact
even-dimensional oriented
 foliation with spin leave and codimension $q$, and $D_F$ be the sub-Dirac operator, then
${\rm lim}_{\varepsilon\rightarrow 0}\varepsilon^{\frac{q}{2}}{\rm
Res}(D^{-n+2}_{F,\varepsilon})$ is proportional to $\int_M[k^F+\Phi(\omega)]d{\rm vol}_{g}$.}\\

\noindent {\bf Remark 1.} When $(M,F)$ is a fibre bundle with
compact fibres, then $\Phi(\omega)=0$, so the
$O(\frac{1}{\varepsilon^{\frac{q}{2}}})$ coefficient of $ {\rm
Res}(D^{-n+2}_{F,\varepsilon})$ is proportional to the sum of the
gravity along fibre. Especially, when $p=0$, Theorem 5.1 is the
classical Kastler-Kalau-Walze theorem. \\

\noindent {\bf Remark 2} By (2.32-2.35) in [LZ], we know that
$$k^{TM,\varepsilon}=k^F+\Phi(\omega)+a_1\varepsilon+a_2\varepsilon^2,\eqno(5.9)$$
where $a_1,a_2$ are functions which are independent of
$\varepsilon.$ So, if $q=2$ (resp. $4$), we will have that ${\rm
lim}_{\varepsilon\rightarrow 0}{\rm Res}(D^{-n+2}_{F,\varepsilon})$
is proportional to $\int_M a_1d{\rm vol}_{g}$ (resp. $\int_M
a_2d{\rm vol}_{g}$).\\

\noindent {\bf Acknowledgement.} The work of the first author was
partially supported by NSF and NSFC. The work of the second author
was supported by Science Foundation for Young Teachers of Northeast
Normal University (No. 20060102). We thank Prof. Weiping Zhang for his generous discussions, especially
the idea of proving a vanishing theorem in Section 2.\\

\noindent{\bf Reference}\\

\noindent [Ac] T. Ackermann, {\it A note on the Wodzicki residue.}
J. Geom. Phys. 20 (1996), no. 4, 404--406.

\noindent [Bi] J. M. Bismut, {\it A local index theorem for
non-K\"{a}hler manifolds.} Math. Ann. 284 (1989), no. 4, 681--699.

 \noindent [Bl] D. Blair, {\it Contact manifolds in Riemannian
geometry}. Lecture Notes in Mathematics, Vol. 509. Springer-Verlag,
Berlin-New York, 1976.

 \noindent [Co] A. Connes {\it Cyclic cohomology and the transverse fundamental class of a foliation.}
 Geometric Methods in Operator Algebras, H. Araki eds., pp. 52-144,
 Pitman Research Notes in Math. Series, vol. 123,1986.

\noindent [FGV] H. Figueroa, J. Gracia-Bond\'{i}a and J.
V\'{a}rilly, {\it Elements of noncommutative geometry},
Birkh\"{a}user Boston, 2001.

\noindent [Ka] D. Kastler, {\it The Dirac operator and gravitation.}
Comm. Math. Phys. 166 (1995), no. 3, 633--643.

\noindent [Ko1] Y. Kordyukov, {\it Noncommutative spectral geometry
of Riemannian foliations.} Manuscripta Math. 94 (1997), no. 1,
45--73.

 \noindent [Ko2] Y. Kordyukov, {\it
Adiabatic limits and spectral geometry of foliations.} Math. Ann.
313 (1999), no. 4, 763--783.

\noindent [KW] W. Kalau, and M. Walze, {\it Gravity, non-commutative
geometry and the Wodzicki residue.} J. Geom. Phys. 16 (1995), no. 4,
327--344.

 \noindent [LK] J. A. L$\acute{o}$pez and  Y. A.
Kordyukov, {\it Adiabatic limits and spectral sequences for
Riemannian foliations}. Geom. Funct. Anal. 10 (2000), no. 5,
977--1027.

\noindent [LM] H. Lawson and M. Michelsohn,{\it Spin geometry.}
Princeton Mathematical Series, 38. Princeton University Press,
Princeton, NJ, 1989.

 \noindent [LMZ] K. Liu, X. Ma and W. Zhang,{\it On
elliptic genera and foliations.} Math. Res. Lett. 8 (2001), no. 3,
361--376.

 \noindent [LZ] K. Liu and W. Zhang, {\it Adiabatic
limits and foliations.} Topology, geometry, and algebra:
interactions and new directions (Stanford, CA, 1999), 195--208,
Contemp. Math., 279, Amer. Math. Soc., Providence, RI, 2001.

\noindent [Pe] R. Petit, {\it ${\rm Spin}\sp c$-structures and Dirac
operators on contact manifolds.} Differential Geom. Appl. 22 (2005),
no. 2, 229--252.

 \noindent [Ru] M. Rumin, {\it Sub-Riemannian limit of the
differential form spectrum of contact manifolds.} Geom. Funct. Anal.
10 (2000), no. 2, 407--452.

\noindent [Wo] M. Wodzicki, {\it Local invariants of spectral
asymmetry.} Invent. Math. 75 (1984), no. 1, 143--177.\\

\indent{ Center of Mathematical Sciences, Zhejiang University
Hangzhou Zhejiang 310027, China and Department of Mathematics,
University of California at Los Angeles, Los Angeles CA 90095-1555,
USA\\} \indent  Email: {\it liu@ucla.edu.cn; liu@cms.zju.edu.cn}\\

 \indent{  School of Mathematics and Statistics,
Northeast Normal University, Changchun Jilin, 130024 and  Center of
Mathematical Sciences, Zhejiang
University, Hangzhou Zhejiang 310027, China }\\
\indent E-mail: {\it wangy581@nenu.edu.cn}\\

\end{document}